\documentclass[11pt]{article}

\usepackage{graphicx, subfigure}
\usepackage[top=1in,bottom=1in,left=1in,right=1in]{geometry}
\usepackage[sort&compress,numbers]{natbib} 
\usepackage{amsfonts, amsmath, amssymb, amsthm, constants}

\usepackage{color}
\definecolor{darkred}{RGB}{100,0,0}
\definecolor{darkgreen}{RGB}{0,100,0}
\definecolor{darkblue}{RGB}{0,0,150}

\usepackage{url}

\newtheorem{thm}{Theorem}
\newtheorem{prp}{Proposition}

\def\beq{\begin{equation}}
\def\eeq{\end{equation}}
\def\beqn{\begin{eqnarray*}}
\def\eeqn{\end{eqnarray*}}
\def\bitem{\begin{itemize}}
\def\eitem{\end{itemize}}
\def\benum{\begin{enumerate}}
\def\eenum{\end{enumerate}}
\def\bmult{\begin{multline*}}
\def\emult{\end{multline*}}
\def\bcenter{\begin{center}}
\def\ecenter{\end{center}}

\newcommand{\thmref}[1]{Theorem~\ref{thm:#1}}
\newcommand{\prpref}[1]{Proposition~\ref{prp:#1}}

\newcommand{\secref}[1]{Section~\ref{sec:#1}}

\def\cN{\mathcal{N}}

\def\cX{\mathcal{X}}

\def\bC{\mathbf{C}}

\def\bI{\mathbf{I}}

\def\bM{\mathbf{M}}

\def\ba{\mathbf{a}}
\def\bb{\mathbf{b}}

\def\bx{\mathbf{x}}

\def\bbP{\mathbb{P}}

\def\bbR{\mathbb{R}}

\newcommand{\E}{\operatorname{\mathbb{E}}}
\renewcommand{\P}{\operatorname{\mathbb{P}}}

\newcommand{\expect}[1]{\mathbb{E}\left(#1\right)}

\newcommand{\var}[1]{\operatorname{Var}\left(#1\right)}

\newcommand{\<}{\langle}
\renewcommand{\>}{\rangle}

\def\eps{\varepsilon}

\begin{document}
\thispagestyle{empty}
% TITLE
\noindent {\sc \LARGE Detecting a Vector Based on Linear Measurements}

% AUTHORS
\bigskip 
\noindent {\Large
Ery Arias-Castro} \\[.05in]
{\large
Department of Mathematics, University of California, San Diego
}

\bigskip
\noindent 
We consider a situation where the state of a system is represented by a real-valued vector $\bx \in \bbR^n$.  Under normal circumstances, the vector $\bx$ is zero, while an event manifests as non-zero entries in $\bx$, possibly few.  Our interest is in designing algorithms that can reliably detect events --- i.e., test whether $\bx = 0$ or $\bx \neq 0$ --- with the least amount of information.  We place ourselves in a situation, now common in the signal processing literature, where information on $\bx$ comes in the form of noisy linear measurements $y = \<\ba, \bx\> + z$, where $\ba \in \bbR^n$ has norm bounded by 1 and $z \in \cN(0,1)$.  We derive information bounds in an active learning setup and exhibit some simple near-optimal algorithms.  
%For example, when $\bx$ has non-negative entries, then a simple algorithm is able to reliably detect with $m$ measurements if
%$
%\sqrt{\frac{m}{n}} \sum_{j=1}^n x_j \to \infty,
%$
%while no algorithm, adaptive or non-adaptive, is able to reliably detect if
%$
%\sqrt{\frac{m}{n}} \sum_{j=1}^n x_j \leq \kappa,
%$
%for an explicit constant $\kappa$.
In particular, our results show that the task of detection within this setting is at once much easier, simpler and different than the tasks of estimation and support recovery.

\medskip


\section{Introduction}
\label{sec:intro}

We consider a situation where the state of a system is represented by a real-valued vector $\bx \in \bbR^n$.  Under normal circumstances, the vector $\bx$ is zero, while an event manifests as non-zero entries in $\bx$, possibly few.  Our interest is in the design of algorithms that reliably detect events --- i.e., test whether $\bx = 0$ or $\bx \neq 0$ --- with the least amount of information.  We assume that we may learn about $\bx$ via noisy linear measurements of the form 
\beq \label{model}
y_i = \<\ba_i, \bx\> + z_i, 
\eeq
where the measurement vectors $\ba_i$'s have Euclidean norm bounded by 1 and the noise $z_i$'s are i.i.d.~standard normal.  Assuming that we may take a limited number of linear measurements, the engineering is in choosing them in order to minimize the false alarm and missed detection rates.  We derive information bounds, establishing some fundamental detection limits relating the signal strength and the number of linear measurements.  The bounds we obtain apply to all adaptive schemes, where we may choose the $i$th measurement vector $\ba_i$ based on the past measurements, i.e., we may choose $\ba_i$ as a function of $(\ba_1, y_1, \dots, \ba_{i-1}, y_{i-1})$.

\subsection{Related work}
Learning as much as possible about a vector based on a few linear measurements is one of the central themes of compressive sensing (CS)~\cite{CRT,OptimalRecovery,Donoho-CS}.  Most of this literature, as it relates to signal processing, has focused on the tasks of estimation and support recovery.  Particularly in surveillance situations, however, it makes sense to perform detection before estimation because, as we shall confirm, reliable detection is possible at much lower signal-to-noise ratios or, equivalently, with much fewer linear measurements than estimation.  This can be achieved with much greater implementation ease and much lower computational cost than standard CS methods based on convex programming.    

The literature on the detection of a high-dimensional signal is centered around the classical normal mean model, based on observations $y_i = x_i + z_i$, where the $z_i$'s are i.i.d.~standard normal.  In this model, only one noisy observation is available per coordinate, so that some assumptions are necessary and the most common one, by far, is that the vector $\bx = (x_1, \dots, x_n)$ is sparse.  This setting has attracted a fair amount attention~\cite{IngsterBook,Ingster99,dj04}, with recent publications allowing adaptive measurements~\cite{haupt2009distilled}.  More recently, a few papers~\cite{ingster2010detection,anova-hc,sparse-detection} extended these results to testing for a sparse coefficient vector in a linear system with the aim of characterizing the detection feasibility.  These papers work with designs having low mutual coherence, for example, assuming that the $\ba_i$'s are i.i.d.~multivariate normal.  As we shall see below, such designs are not always desirable.
We also mention~\cite{CS-detection}, which assumes that an estimator $\widehat{\bx}$ of $\bx$ is available and examines the performance of the test based on $\<\widehat{\bx}, \bx\>$; and~\cite{meng-sparse}, which proposes a Bayesian approach for the detection of sparse signals in a sensor network for which the design matrix is assumed to have some polynomial decay in terms of the distance between sensors.

We mention that the present paper may be seen as a companion paper to~\cite{adaptiveCS} which considers the tasks of estimation and support recovery in the same setting.
   
%We emphasize at this point that our results apply equally to sparse and non-sparse settings.

\subsection{Notation and terminology}

Our detection problem translates into a hypothesis testing problem $H_0: \bx = 0$ versus $H_1: \bx \in \cX$, for some subset $\cX \subset \bbR^n \setminus \{0\}$.  A test procedure based on $m$ measurements of the form \eqref{model} is a binary function of the data, i.e., $T = T(\ba_1, y_1, \dots, \ba_{m}, y_{m})$, with $T = \eps \in \{0,1\}$ indicating that $T$ favors $H_\eps$.  The (worst-case) risk of a test $T$ is defined as 
$$
\gamma(T) := \bbP_0(T = 1) + \max_{\bx \in \cX} \,  \bbP_\bx(T = 0),
$$
where $\bbP_\bx$ denotes the distribution of the data when $\bx$ is the true underlying vector.  With a prior $\pi$ on the set of alternatives $\cX$, the corresponding average (Bayes) risk is defined as 
$$
\gamma_\pi(T) := \bbP_0(T = 1) + \E_\pi \bbP_\bx(T = 0),
$$
where $\E_\pi$ denotes the expectation under $\pi$.  Note that for any prior $\pi$ and any test procedure $T$,
\beq \label{risks}
\gamma(T) \geq \gamma_\pi(T).
\eeq 
For a vector $\ba = (a_1, \dots, a_k)$, 
\[
\|\ba\| = \big( \sum_j a_j^2 \big)^{1/2}, \qquad |\ba| = \sum_j |a_j|,
\]
and $\ba^T$ denote its transpose.  For a matrix $\bM$, 
\[
\|\bM\|_{\rm op} = \sup_{\ba \neq 0} \frac{\|\bM \ba\|}{\|\ba\|}.
\]
Everywhere in the paper, $\bx = (x_1, \dots, x_n)$ denotes the unknown vector, while ${\bf 1}$ denotes the vector with all coordinates equal to 1 and dimension implicitly given by the context.

\subsection{Content}

In \secref{nonneg} we focus on vectors $\bx$ with non-negative coordinates.  This situation leads to an exceedingly simple, yet near-optimal procedure based on a measurement scheme that is completely at odds with what is commonly used in CS.  In \secref{general} we treat the case of a general vector $\bx$ and derive another simple, near-optimal procedure.  In both cases, the methods we suggest are non-adaptive --- in the sense that the measurement vectors are chosen independently of the observations --- yet perform nearly as well as any adaptive method.  In \secref{discussion} we discuss our results and important extensions, particularly to the case of structured signals.

\section{Vectors with non-negative entries}
\label{sec:nonneg}

Vectors with non-negative entries may be relevant in image processing, for example, where the object to be detected is darker (or lighter) than the background.  As we shall see, detecting such a vector is essentially straightforward in every respect.  In particular, the use of low-coherence designs is counter-productive in this situation.

The first thing that comes to mind, perhaps, is gathering strength across coordinates by measuring $\bx$ with the constant vector ${\bf 1}/\sqrt{n}$.  And, with a budget of $m$ measurements, we simply take this measurement $m$ times. 
%\beq \label{nonneg}
%y_i = \frac1{\sqrt{n}} |\bx| + z_i, \ i = 1, \dots, m.
%\eeq  

\begin{prp} \label{prp:nonneg-ub}
Suppose we take $m$ measurements of the form \eqref{model} with $\ba_i = {\bf 1}/\sqrt{n}$ for all $i$.  Consider then the test that rejects when 
$$
\sum_{i=1}^m y_i > \tau \sqrt{m},
$$
where $\tau$ is some critical value.  Its risk against a vector $\bx$ is equal to
$$
1 - \Phi(\tau) + \Phi(\tau - \sqrt{m/n} |\bx|),
$$
where $\Phi$ is the standard normal distribution function.  In particular, if $\tau = \tau_n \to \infty$, this test has vanishing risk against alternatives satisfying $\sqrt{m/n} |\bx| - \tau_n \to \infty$.
\end{prp}

Since we may chose $\tau_m \to \infty$ as slowly as we wish, in essence, the simple sum test based on repeated measurements from the constant vector has vanishing risk against alternatives satisfying $\sqrt{m/n} |\bx| \to \infty$.

\begin{proof}
The result is a simple consequence of the fact that
$$
\frac1{\sqrt{m}} \sum_{i=1}^m y_i \sim \cN(\sqrt{m/n} |\bx|, 1).
$$
\end{proof}

Although the choice of measurement vectors and the test itself are both exceedingly simple, the resulting procedure comes close to achieving the best possible performance in this particular setting, as the following information bound reveals.

\begin{thm} \label{thm:nonneg-lb}
Let $\cX(\mu, S)$ denote the set of vectors in $\bbR^n$ having exactly $S$ non-zero entries all equal to $\mu > 0$.  Based on $m$ measurements of the form \eqref{model}, possibly adaptive, any test for $H_0: \bx = 0$ versus $H_1: \bx \in \cX(\mu, S)$ has risk at least $1 - \sqrt{m/(8n)} S \mu$. 
\end{thm}

In particular, the risk against alternatives $\bx \in \cX(\mu, S)$ with $\sqrt{m/n} |\bx| = \sqrt{m/n} S \mu \to 0$, goes to 1 uniformly over all procedures.

\begin{proof}
The standard approach to deriving uniform lower bounds on the risk is to put a prior on the set of alternatives and use \eqref{risks}.  We simply choose the uniform prior on $\cX(\mu, S)$, which we denote by $\pi$.  The hypothesis testing problem reduces to $H_0: \bx = 0$ versus $H_1: \bx \sim \pi$, for which the likelihood ratio test is optimal by the Neyman-Pearson fundamental lemma.  The likelihood ratio is defined 
$$
L := \frac{\bbP_\pi(\ba_1, y_1, \dots, \ba_{m}, y_{m})}{\bbP_0(\ba_1, y_1, \dots, \ba_{m}, y_{m})} = \E_\pi \exp\left(\sum_{i = 1}^m y_i (\ba_i^T \bx) - (\ba_i^T \bx)^2/2\right),
$$
where $\E_\pi$ denotes the expectation with respect to $\pi$, and the related test is $T = \{L > 1\}$.  It has risk equal to 
\beq \label{gamma}
\gamma_\pi(T) = 1 - \frac12 \|\bbP_\pi - \bbP_0\|_{\rm TV},
\eeq
where $\P_\pi := \E_\pi \P_\bx$ --- the $\pi$-mixture of $\P_\bx$ --- and $\| \cdot \|_{\rm TV}$ is the total variation distance.  
By Pinsker's inequality
\beq \label{pinsker}
\|\bbP_\pi - \bbP_0\|_{\rm TV} \leq \sqrt{K(\bbP_0, \bbP_\pi)/2},
\eeq
where $K(\bbP_0, \bbP_\pi)$ denotes the Kullback-Leibler divergence.  
We have 
\begin{eqnarray}
K(\bbP_0, \bbP_\pi) &=& - \E_0 \log L \\
&\leq& \E_\pi \sum_{i = 1}^m  \E_0 \left(y_i (\ba_i^T \bx) - (\ba_i^T \bx)^2/2\right) \\ 
&=& \E_\pi \sum_{i = 1}^m \E_0 (\ba_i^T \bx)^2/2 \\
&=& \sum_{i = 1}^m \E_0 \big(\ba_i^T \bC \ba_i \big) \label{sum1} \\
&\le& m \|\bC\|_{\rm op}, \label{sum2}
\end{eqnarray}
where $\bC = (c_{j k}) := \E_\pi (\bx \bx^T)$.  The first line is by definition; the second is by definition of $\P_\bx/\P_0$, by the application of Jensen's inequality justified by the convexity of $x \to - \log x$, and by Fubini's theorem; the third is by independence of $\ba_i$, $y_i$ and $\bx$ (under $\P_0$), and by the fact that $\E(y_i) = 0$; the fourth is by independence of $\ba_i$ and $\bx$ (under $\P_0$) and by Fubini's theorem; the fifth is because $\|\ba_i\| \le 1$ for all $i$.

Since under $\pi$ the support of $\bx$ is chosen uniformly at random among subsets of size $S$, we have 
\[c_{j j} = \mu^2 \P_\pi(x_j \neq 0) = \mu^2 \cdot \frac{S}n, \quad \forall j,\]
and
\[c_{j k} = \mu^2 \P_\pi(x_j \neq 0, x_k \neq 0) = \mu^2 \cdot \frac{S}n \cdot \frac{S-1}{n-1}, \qquad j \neq k.\]
This simple matrix has operator norm $\|\bC\|_{\rm op} = \mu^2 S^2/n$.

Coming back to the divergence, we therefore have
\[K(\bbP_0, \bbP_\pi) \le m \cdot \mu^2 S^2/n,\] 
and returning to \eqref{gamma} via \eqref{pinsker}, we bound the risk of the likelihood ratio test as follows
$$
\gamma(T)  \geq 1 - \sqrt{K(\bbP_0, \bbP_\pi)/8} \geq 1 - \sqrt{m/(8n)} S \mu.
$$
\end{proof}

With \prpref{nonneg-ub} and \thmref{nonneg-lb}, we conclude that the following is true in a minimax sense:
\begin{quote}
{\it Reliable detection of a nonnegative vector $\bx \in \bbR^n$ from $m$ noisy linear measurements is possible if $\sqrt{m/n} |\bx| \to \infty$ and impossible if $\sqrt{m/n} |\bx| \to 0$.}  
\end{quote}

\section{General vectors}
\label{sec:general}

When dealing with arbitrary vectors, the measurement vector ${\bf 1}/\sqrt{n}$ may not be appropriate.  In fact, the resulting procedure is completely insensitive to vectors $\bx$ such that $\<{\bf 1}, \bx\> = 0$.  Nevertheless, if one selects a measurement vector $\ba$ from the Bernoulli ensemble --- i.e., with independent entries taking values $\pm 1/\sqrt{n}$ with equal probability --- then on average, $\<\ba, \bx\>$ is of the order of $\|\bx\|/\sqrt{n}$.  This is true when the number of non-zero entries in $\bx$ grows with the dimension $n$; if we repeat the process a few times, it becomes true for any fixed vector $\bx$.  

\begin{prp} \label{prp:general-ub}
Sample $\bb_1, \dots, \bb_{h_m}$ independently from the Bernoulli ensemble, with $h_m \to \infty$ slowly, and take $m$ measurements of the form \eqref{model} with $\ba_i = \bb_{s}$ for $i \in I_s := [(m/h_m)(s-1)+1, (m/h_m)s), \ s = 1, \dots, h_m$.  Consider then the test that rejects when 
\beq \label{general-ub}
\sum_{s=1}^{h_m} \left(\sum_{i \in I_s} y_i\right)^2 > m (1 + \tau_m/\sqrt{h_m}),
\eeq
where $\tau_m \to \infty$.  When $m \to \infty$, its risk against a vector $\bx$ --- averaged over the Bernoulli ensemble --- vanishes if $(m/n) \|\bx\|^2 \geq 2\tau_m \sqrt{h_m}$.
\end{prp}

Since we may take $h_m$ and $\tau_m$ increasing as slowly as we please, in essence, the test is  reliable when $(m/n) \|\bx\|^2 \to \infty$.  Compared with repeatedly measuring with the constant vector ${\bf 1}/\sqrt{n}$ as studied in \prpref{nonneg-ub}, there is a substantial loss in power when $|\bx|^2$ is much larger than $\|\bx\|^2$.  For example, when $\bx$ has $S$ non-zero entries al equal to $\mu > 0$, $|\bx|^2 = S \|\bx\|^2$.

\begin{proof}
For simplicity, assume that $m/h_m$ is an integer and fix $\bx$ throughout.  For short, let 
$$
Y_s = \sum_{i \in I_s} y_i = (m/h_m) \<\bb_s, \bx\> + \sqrt{m/h_m} Z_s, \quad Z_s := \sqrt{h_m/m}\sum_{i \in I_s} z_i.
$$
Note that the $Z_s$'s are i.i.d.$\sim \cN(0, 1)$, while the $\<\bb_s, \bx\>$'s are i.i.d.~with mean zero, variance $\|\bx\|^2/n$ and fourth moment bounded by $6 \|\bx\|^4/n^2$ --- which is immediate using the fact that the coordinates of $\bb_s$ are i.i.d.~taking values $\pm 1/\sqrt{n}$ with equal probability.  
Proceeding in an elementary way, we have
$$
\expect{\sum_{s=1}^{h_m} Y_s^2} = (m^2/h_m) \expect{\<\bb_1, \bx\>^2} +  m \expect{Z_1^2} = (m^2/h_m) \|\bx\|^2/n + m,
$$
and 
\beqn
\var{\sum_{s=1}^{h_m} Y_s^2} &=& (m^4/h_m^3) \expect{\<\bb_1, \bx\>^4} + 6 (m^3/h_m^2) \expect{\<\bb_1, \bx\>^2} \expect{Z_1^2} + (m^2/h_m) \expect{Z_1^4} \\
&=& 6 (m^4/h_m^3) \|\bx\|^4/n^2 + 6 (m^3/h_m^2) \|\bx\|^2/n + 3 (m^2/h_m).
\eeqn
Therefore, by Chebyshev's inequality, the probability of \eqref{general-ub} under the null is bounded from above by $3/\tau_m^2 \to 0$.  Similarly, the probability of \eqref{general-ub} {\em not} happening under an alternative $\bx$ satisfying $(m/n) \|\bx\|^2 \geq 2\tau_m \sqrt{h_m}$ is bounded from above by
$$
\frac{6 (m^4/h_m^3) \|\bx\|^4/n^2 + 6 (m^3/h_m^2) \|\bx\|^2/n + 3 (m^2/h_m)}{((m^2/h_m) \|\bx\|^2/n - m \tau_m/\sqrt{h_m})^2} \leq  \frac{24}{h_m} + \frac{24}{\tau_m \sqrt{h_m}} + \frac{3}{\tau_m^2} \to 0.
$$ 
\end{proof}

Again, this relatively simple procedure nearly achieves the best possible performance.

\begin{thm} \label{thm:general-lb}
Let $\cX^\pm(\mu, S)$ denote the set of vectors in $\bbR^n$ having exactly $S$ non-zero entries all equal to $\pm \mu$.  Based on $m$ measurements of the form \eqref{model}, possibly adaptive, any test for $H_0: \bx = 0$ versus $H_1: \bx \in \cX^\pm(\mu, S)$ has risk at least $1 - \sqrt{S m/(8n)} \mu$. 
\end{thm}

In particular, the risk against alternatives $\bx \in \cX^\pm(\mu, S)$ with $(m/n) \|\bx\|^2 = (m/n) S \mu^2 \to 0$, goes to 1 uniformly over all procedures.

\begin{proof}
Again, we choose the uniform prior on $\cX^\pm(\mu, S)$.  The proof is then completely parallel to that of \thmref{nonneg-lb}, now with $\bC = \mu^2 (S/n) \bI$ --- since the signs of the nonzero entries of $\bx$ are i.i.d.~Rademacher --- so that $\|\bC\|_{\rm op} = \mu^2 S/n$.
\end{proof}

With \prpref{general-ub} and \thmref{general-lb}, we conclude that the following is true in a minimax sense:
\begin{quote}
{\it Reliable detection of a vector $\bx \in \bbR^n$ from $m$ noisy linear measurements is possible if $\sqrt{m/n} \|\bx\| \to \infty$ and impossible if $\sqrt{m/n} \|\bx\| \to 0$.}  
\end{quote}

\section{Discussion}
\label{sec:discussion}

In this short paper, we tried to convey some very basic principles about detecting a high-dimensional vector with as few linear measurements as possible.  First, when the vector has non-negative entries, repeatedly sampling from the constant vector ${\bf 1}/\sqrt{n}$ is near-optimal.  Second, when the vector is general but sparse, repeatedly sampling from a few measuring vectors drawn from a standard random (e.g., Bernoulli) ensemble is also near-optimal.  In both cases, choosing the measuring vectors adaptively does not bring a substantial improvement.  And, moreover, sparsity does not help, in the sense that the detection rates depend on $|\bx|$  and $\|\bx\|$, respectively.

\subsection{A more general adaptive scheme}
\label{sec:energy}

Suppose we may take as many linear measurements of the form \eqref{model} as we please (possibly an infinite number), with the only constraint being on the total measurement energy 
\beq \label{energy}
\sum_i \|\ba_i\|^2 \le m.
\eeq
(Note that $m$ is no longer constrained to be an integer.)  This is essentially the setting considered in~\cite{haupt2009distilled,haupt-compressive}, and clearly, the setup we studied in the previous sections satisfies this condition.  
So what can we achieve with this additional flexibility?  

In fact, the same results apply.  The lower bounds in \thmref{nonneg-lb} and \thmref{general-lb} are proved in exactly the same way.  (We effectively use \eqref{energy} to go from \eqref{sum1} to \eqref{sum2}, and this is the only place where the constraints on the number and norm of the measurement vectors are used.)  Of course, \prpref{nonneg-ub} and \prpref{general-ub} apply since the measurement schemes used there satisfy \eqref{energy}.  However, in this special case they could be simplified.  For instance, in \prpref{nonneg-ub} we could take one measurement with the constant vector $\sqrt{m/n} \, {\bf 1}$.

\subsection{Detecting structured signals}

The results we derived are tailored to the case where $\bx$ has no known structure.  What if we know a priori that the signal $\bx$ has some given structure?  The most emblematic case is when the support of $\bx$ is an interval of length $S$.  In the classical setting where each coordinate of $\bx$ is observed once, the scan statistic (aka generalized likelihood ratio test) is the tool of choice~\cite{MGD}.  How does the story change in the setting where adaptive linear measurements in the form of \eqref{model} can be taken?  

Perhaps surprisingly, knowing that $\bx$ has such a specific structure does not help much.  Indeed, \thmref{nonneg-lb} and \thmref{general-lb} are proved in the same way.  In the case of non-negative vectors, we use the uniform prior on vectors with support an interval of length $S$ and nonzero entries all equal to $\mu$, and the proof is identical, except for the matrix $\bC$, which now has coefficients $c_{jk} = \mu^2 \max(S - |j-k|,0)/n$ for all $j,k$.  Because $\bC$ is symmetric, we have
\beq
\|\bC\|_{\rm op} \le \max_j \sum_k |c_{jk}| = \mu^2 S^2/n,
\eeq
which is exactly the same bound as before.  In the general case, the arguments are really identical, except that we use uniform prior on vectors with support an interval of length $S$ and nonzero entries all equal to $\mu$ in absolute value.  (Here the matrix $\bC$ is exactly the same.)  Of course, \prpref{nonneg-ub} and \prpref{general-ub} apply here too, so the conclusions are the same.  Here too, these conclusions hold in the more general setup with measurements satisfying \eqref{energy}.

To appreciate how powerful the ability to take linear measurements in the form of \eqref{model} with the constraint \eqref{energy} really is, let us stay with the same task of detecting an interval of length $S$ with a positive mean.  On the one hand, we have the simple test based on $\sum_i y_i$ studied in \prpref{nonneg-ub}.  On the other hand, we have the scan statistic 
\[
\max_{t} \sum_{i = t}^{t+S-1} y_i,
\] 
with observations of the form  
\beq \label{classic_model}
y_i = x_i + \sigma z_i, \quad \sigma := \sqrt{n/m}.
\eeq
While the former requires $\sqrt{m/n} |\bx| \to \infty$ to be asymptotically powerful, the scan statistic requires 
\[\varliminf \sqrt{m/n} |\bx| \cdot (S \log^+(n/S))^{-1/2} \ge \sqrt{2},\]  
where $\log^+(x) := \max(\log x, 1)$.  With observations provided in the form of \eqref{classic_model}, this is asymptotically optimal~\cite{MGD}.  Note that \eqref{classic_model} is a special case of \eqref{energy}.  Hence, the ability of taking measurements of the form \eqref{model} allows to detect structured signals that are potentially much weaker, without a priori knowledge of the structure and with much simpler algorithms.  Hardware that is able to take linear measurements such as \eqref{model} is currently being developed~\cite{duarte2008single}.

\subsection{A comparison with estimation and support recovery}

The results we obtain for detection are in sharp contrast with the corresponding results in estimation and support recovery.  Though, by definition, detection is always easier, in most other settings it is not that much easier.  For example, take the normal mean model described in the Introduction, assuming $\bx$ is sparse with $S$ coefficients equal to $\mu > 0$.  In the regime where $S = p^{1-\beta}, \ \beta \in (1/2,1)$, detection is impossible when $\mu \leq \sqrt{2 r \log n}$ with $r < \rho_1(\beta)$, while support recovery is possible when $\mu \geq \sqrt{2 r \log n}$ with $r > \rho_2(\beta)$, for a fixed functions $\rho_1, \rho_2: (1/2,1) \to (0,\infty)$~\cite{IngsterBook,Ingster99,dj04}.  So the difference is a constant factor in the per-coordinate amplitude.  In the setting we consider here, we are able to detect at a much smaller signal-to-noise ratio than what is required for estimation or support recovery, which nominally require at least $m \geq S$ measurements regardless of the signal amplitude, where $S$ is the number of nonzero entries in $\bx$.  In fact, \cite{adaptiveCS} shows that reliable support recovery is impossible unless $\mu$ is of order at least $\sqrt{n/m}$.  In detection, however, we saw that $m = 1$ measurement may suffice if the signal amplitude is large enough, which can be smaller than $\sqrt{n/m}$ by a factor of $S$ or $\sqrt{S}$ in the nonnegative and general cases respectively.  Therefore, having the ability to take linear measurements of the form \eqref{model} in a surveillance setting, it makes sense to perform detection as described here before estimation (identification) or support recovery (localization) of the signal.

\subsection{Possible improvements}

Though we provided simple algorithms that nearly match information bounds, there might be room for improvement.  For one thing, it might be possible to reliably detect when, say, $\sqrt{m/n} |\bx|$ is sufficiently large --- for the case where $x_j \geq 0$ for all $j$ --- without necessarily tending to infinity.  A good candidate for this might be the Bayesian algorithm proposed in~\cite{4518814}.

More importantly, in the general case of \secref{general}, we might want to design an algorithm that detects any fixed $\bx$ with high-probability, without averaging over the measurement design.  This averaging may be interpreted in at least two ways: 
\renewcommand{\theenumi}{(A\arabic{enumi})}
\renewcommand{\labelenumi}{\theenumi}
\begin{enumerate}
\item \label{a1} If we were to repeat the experiment many times, each time choosing new measurement vectors and corrupting the measurements with new noise, then for a fixed vector $\bx$, in most instances the test would be accurate.
\item \label{a2} Given the amplitudes $|x_j|, \ j = 1, \dots, n$, for most sign configurations the test will be accurate.
\end{enumerate}
Interpretation \ref{a1} is controversial as we do not repeat the experiment, which would amount to taking more samples.  And interpretation \ref{a2} raises the issue of robustness to any sign configuration.  
One way --- and the only way we know of --- to ensure this robustness is to use a CS-like sampling scheme, i.e., choosing $\ba_1, \dots, \ba_m$ in \eqref{model} such that the matrix with these rows satisfies RIP-like properties.  This setting is studied in detail in~\cite{anova-hc}, which in a nutshell says the following.  Take measurement vectors from the Bernoulli ensemble, say, but hold the measurement design fixed.  This is just a way to build a measurement matrix satisfying the RIP and with low mutual coherence.  In particular, this  requires that $m$ is of order at least $S \log n$, though what follows assumes that $m \gg S (\log n)^3$.  Based on such measurements, the test based on $\sum_i y_i^2$ is able to detect when $(\sqrt{m}/n) \|\bx\|^2 \to \infty$, which is more stringent than what is required in \prpref{general-ub}; while the test based on $\max_{j=1, \dots, n} |\sum_i a_{ij} \, y_i|$ is able to detect when $\liminf \sqrt{m/n} \max_j |x_j| (\log n)^{-1/2} > \sqrt{2}$, which, except for the log factor, is what is required for support recovery.  And this is essential optimal, as shown in~\cite{anova-hc}.

%\vspace{.5in}
%\begin{center}
%{\Large Acknowledgements}
%\end{center}

\subsection*{Acknowledgements}

The author would like to thank Emmanuel Cand\`es for stimulating discussions and Rui Castro for clarifying the set up described in \secref{energy}.  This work was partially supported by a grant from the Office of Naval Research (N00014-09-1-0258). 

% bibliograhy
%\renewcommand{\refname}{\begin{center}{\normalfont \huge References}\end{center}}
\bibliographystyle{abbrv}
\bibliography{ref}

\end{document}